\documentclass{amsart}
\usepackage[T1]{fontenc}
\usepackage[utf8]{inputenc}
\usepackage[english]{babel}
\usepackage{csquotes}

\usepackage{biblatex}
\renewbibmacro{in:}{%
        \ifentrytype{article}{}{\printtext{\bibstring{in}\intitlepunct}}}

\usepackage[hang,flushmargin]{footmisc} 

\usepackage{amsmath}
\usepackage{amssymb}
\usepackage{amsthm}

\usepackage{abstract}
\addto\captionsenglish{}

\newtheorem{defn}{Definition}[section]
\newtheorem{theo}[defn]{Theorem}
\newtheorem{prop}[defn]{Proposition}
\newtheorem{lemma}[defn]{Lemma}

\newtheorem*{quest}{Question}

\newcommand{\Q}{\mathbb{Q}}
\newcommand{\R}{\mathbb{R}}

\newcommand{\func}[3]{#1:#2\rightarrow#3}

\renewcommand{\epsilon}{\varepsilon}
\renewcommand{\theta}{\vartheta}

\title[Generic unif.\ cont.\ mappings on unbounded hyperbolic spaces]{Generic uniformly continuous mappings on unbounded hyperbolic spaces}
\author{Davide Ravasini}
\address{Mathematisches Institut, Universit\"{a}t Leipzig \newline\indent Augustusplatz 10, 04109 Leipzig, Germany}
\email{davide.ravasini@uni-leipzig.de}

\begin{document}
\maketitle
\let\thefootnote\relax\footnote{\today \newline \indent\emph{2020 Mathematics Subject Classification:} 54E50, 54E52. \newline
}

\begin{abstract}
\noindent \textsc{Abstract}. We consider a complete, unbounded hyperbolic metric space $X$ and a concave, nonzero and nondecreasing function $\omega:[0,+\infty)\to[0,+\infty)$ with $\omega(0)=0$ and study the space $\mathcal{C}_\omega(X)$ of uniformly continous self-mappings on $X$ whose modulus of continuity is bounded from above by $\omega$. We endow $\mathcal{C}_\omega(X)$ with the topology of uniform convergence on bounded sets and prove that the modulus of continuity of a generic mapping in $\mathcal{C}_\omega(X)$, in the sense of Baire categories, is precisely $\omega$. Some related results in spaces of bounded mappings and in the topology of pointwise convergence are also discussed. This note can be seen as a completion of various results due to F.\ Strobin, S.\ Reich, A.J.\ Zaslavski, C.\ Bargetz and D.\ Thimm.
\end{abstract}

\section{Introduction}
The starting point of our discussion is the set of all nonexpansive mappings $\func{f}{C}{C}$, where $C$ is a closed and convex set in a Banach space $X$. These are mappings $f$ such that ${\|f(x)-f(y)\|}_X\leq{\|x-y\|}_X$ for every $x,y\in C$. In \cite{dbm_1989}, F.S.\ De Blasi and J.\ Myjak showed that, if $X$ is a Hilbert space and $C$ is bounded, the set of strict contractions $\func{f}{C}{C}$ is a $\sigma$-porous subset of the space of all nonexpansive mappings on $C$, endowed with the topology of uniform convergence. Topologically, this means that the size of the set of strict contractions is negligible and that most nonexpansive mappings are not in fact strict contractions. Although Banach's classical fixed point theorem can only be applied to strict contractions, De Blasi and Myjak also showed that the set of nonexpansive mappings on $C$ without a fixed point is $\sigma$-porous whenever $C$ is a closed, convex, bounded set in a Banach space.

S.\ Reich and A.J.\ Zaslavski later showed in \cite{rz_2001} that if $C$ is a closed, convex, bounded subset of a Banach space, then, with the exception of a $\sigma$-porous set, every nonexpansive mapping on $C$ is contractive in the sense of Rakotch (see \cite{rak_1962}). Rakotch contractions do have unique fixed points, therefore the result of Reich and Zaslavski generalises the work of De Blasi and Myjak. On the other hand, C.\ Bargetz and M.\ Dymond showed that the set of strict contractions is $\sigma$-porous also in the case of a general Banach space (\cite{cm_2016}). This result was later extended by C.\ Bargetz, M.\ Dymond and S.\ Reich in \cite{bdr_2017} to a larger class of geodesic metric spaces.

Motivated by the work of De Blasi and Myjak, in \cite{fs_2012} F.\ Strobin investigates more general spaces of continuous Hilbert space mappings. It follows as a consequence of his results that the situation on unbounded domains is radically different. Given a closed, convex and unbounded set $K$ in a Hilbert space and a concave, nonzero and nondecreasing function $\func{\omega}{[0,+\infty)}{[0,+\infty)}$ with $\omega(0)=0$, it is possible to define the space $\mathcal{C}_\omega(K)$ of uniformly continuous mappings $\func{f}{K}{K}$ with $\omega_f\leq\omega$, where $\omega_f$ is the modulus of continuity of $f$, and to endow $\mathcal{C}_\omega(K)$ with the topology of uniform convergence on bounded sets. Strobin showed that, with the exception of a meagre set, every mapping $f\in\mathcal{C}_\omega(K)$ is such that $\omega_f=\omega$. In the case of nonexpansive mappings, i.e.\ when $\omega(s)=s$ for every $s$, this result implies that Rakotch contractions form a meagre subset of the space of all nonexpansive mappings. It was then shown in \cite{brt_2023} by C.\ Bargetz, S.\ Reich and D.\ Thimm that Rakotch contractions form a meagre subset of the space of nonexpansive mappings also in the more general class of complete, unbounded hyperbolic metric spaces. However, \cite{brt_2023} does not generalise Strobin's result for uniformly continous mappings to hyperbolic spaces, as it only considers nonexpansive mappings. The proofs in \cite{dbm_1989} and \cite{fs_2012} heavily rely on the Kirszbraun-Valentine Theorem and the Grünbaum-Zarantonello Theorem on the extension of Lipschitz and uniformly continuous mappings respectively, which are available exclusively in Hilbert spaces. Therefore, any possible generalisation of these results to more general spaces usually requires rather different approaches.

The aim of the present work is to prove a generalisation of Strobin's theorem to the case of a complete, unbounded hyperbolic space $X$. Given a concave, nonzero, nondecreasing function $\func{\omega}{[0,+\infty)}{[0,+\infty)}$, we show that a generic mapping $f\in\mathcal{C}_\omega(X)$ is such that $\omega_f=\omega$. We also show that in spaces of bounded, uniformly continuous mappings, endowed with the topology of uniform convergence, the situation is completely different, in that the generic mapping $f$ satisfies the inequality $\omega_f(s)<\omega(s)$ at every point $s\in(0,+\infty)$. Finally, we discuss the validity of the main theorem with respect to the topology of pointwise convergence under the assumption that the underlying space is separable and present some open questions.

\section{Setting and notation}
\label{sec:setting}
For a more detailed introduction to hyperbolic spaces, see \cite{brt_2023}. We report here the most important tools and facts we need. The notation we use is standard for metric spaces. In particular, if $(X,d)$ is a metric space, $x\in X$ and $r>0$, we denote by $B(x,r)$ the open ball of radius $r$ centred around $x$.

A metric space $(X,\rho)$ is \emph{geodesic} if for every $x,y\in X$ there is an isometry $\func{\gamma}{[0,\rho(x,y)]}{X}$ with $\gamma(0)=x$ and $\gamma(\rho(x,y))=y$. We call $\gamma$ a \emph{geodesic} from $x$ to $y$. If, for every $x,y\in X$, the geodesic from $x$ to $y$ is unique, we say that $X$ is \emph{uniquely geodesic}. In general, we do not require our spaces to be uniquely geodesic, as this definition is not even satisfied by Banach spaces, unless they are strictly convex. What we require is a \emph{coherent system} of geodesics on $X$. This is a family $\Gamma$ of geodesics on $X$ with the following properties.
\begin{enumerate}
\item For every $x,y\in X$, there is a unique geodesic $\gamma\in\Gamma$ from $x$ to $y$.
\item If $\gamma\in\Gamma$ is the unique geodesic from $x$ to $y$, then 
\[ \func{\overline{\gamma}}{[0,\rho(x,y)]}{X},\quad\overline{\gamma}(t)=\gamma\bigl(\rho(x,y)-t\bigr) \]
is the unique geodesic in $\Gamma$ from $y$ to $x$
\item If $\gamma\in\Gamma$ is the unique geodesic from $x$ to $y$, then for any $t_1,t_2\in[0,\rho(x,y)]$ with $t_2\geq t_1$ the function
\[ \func{\delta}{[0,t_2-t_1]}{X},\quad\delta(t)=\gamma(t+t_1) \] 
is the unique geodesic in $\Gamma$ from $z_1=\gamma(t_1)$ to $z_2=\gamma(t_2)$.
\end{enumerate}

If $\Gamma$ is a coherent system of geodesics on $X$, for any $x,y\in X$ we define the \emph{geodesic segment} between $x$ and $y$ as $[x,y]=\gamma([0,\rho(x,y)])$, where $\gamma\in\Gamma$ is the unique geodesic from $x$ to $y$. Properties (2) and (3) in the definition of $\Gamma$ clearly imply that $[z,w]\subseteq[x,y]$ for every $x,y,z,w\in X$ with $z,w\in[x,y]$. Conversely, if $\Gamma$ satisfies (1) and $[z,w]\subseteq[x,y]$ for every $x,y,z,w\in X$ with $z,w\in[x,y]$, then $\Gamma$ must satisfy (2) and (3). Indeed, if $\gamma,\overline{\gamma}\in\Gamma$ are the unique geodesics from $x$ to $y$ and from $y$ to $x$ respectively, then by $[x,y]=[y,x]$ we have that $\func{\gamma^{-1}\circ\overline{\gamma}}{[0,\rho(x,y)]}{[0,\rho(x,y)]}$ is an isometry with $\gamma^{-1}\circ\overline{\gamma}(0)=\rho(x,y)$ and $\gamma^{-1}\circ\overline{\gamma}(\rho(x,y))=0$. This holds only if $\gamma^{-1}\circ\overline{\gamma}(t)=\rho(x,y)-t$ for every $t\in[0,\rho(x,y)]$. Similarly, if $\gamma\in\Gamma$ is the unique geodesics from $x$ to $y$, $t_1,t_2\in[0,\rho(x,y)]$ with $t_2\geq t_1$ and $\delta\in\Gamma$ is the unique geodesics from $z_1=\gamma(t_1)$ to $z_2=\gamma(t_2)$, then by $[z_1,z_2]\subseteq[x,y]$ we have that $\func{\gamma^{-1}\circ\delta}{[0,t_2-t_1]}{[0,\rho(x,y)]}$ is an isometry with $\gamma^{-1}\circ\delta(0)=t_1$ and $\gamma^{-1}\circ\delta(t_2-t_1)=t_2$. This holds only if $\gamma^{-1}\circ\delta(t)=t_1+t$ for every $t\in[0,t_2-t_1]$. 

Having a coherent system of geodesics allows us to define convex combinations. For every $x,y\in X$ and every $t\in[0,1]$ we set
\[ (1-t)x\oplus ty=\gamma\bigl(t\rho(x,y)\bigr), \]
where $\gamma\in\Gamma$ is the unique geodesic from $x$ to $y$. We further say that a set $C\subseteq X$ is \emph{convex} if $[x,y]\subseteq C$ for every $x,y\in C$. Note the following remarks, which will be used multiple times throughout the proofs.
\begin{lemma}
\label{lemma:segment}
Let $(X,\rho)$ be a metric space with a coherent system of geodesics $\Gamma$. Then the following assertions hold.
\begin{enumerate}
\item For every $x,y\in X$ and every $\lambda_1,\lambda_2\in[0,1]$ we have 
\[ \rho\bigl((1-\lambda_1)x\oplus\lambda_1 y,(1-\lambda_2)x\oplus\lambda_2 y\bigr)=|\lambda_1-\lambda_2|\rho(x,y). \]
\item For every $x,y\in X$ and every $z\in[x,y]$ we have $\rho(x,y)=\rho(x,z)+\rho(z,y)$.
\end{enumerate}
\end{lemma}
\begin{proof}
To show (1), let $\gamma\in\Gamma$ be the unique geodesics from $x$ to $y$. Since $\gamma$ is an isometry we have 
\begin{align*} 
\rho\bigl((1-\lambda_1)x\oplus\lambda_1 y,(1-\lambda_2)x\oplus\lambda_2 y\bigr) &= \rho\bigl(\gamma(\lambda_1\rho(x,y)),\gamma(\lambda_2\rho(x,y))\bigr)= \\
										&= |\lambda_1\rho(x,y)-\lambda_2\rho(x,y)|= \\
										&= |\lambda_1-\lambda_2|\rho(x,y).
\end{align*}

To show (2), let $\lambda\in[0,1]$ be such that $z=(1-\lambda)x\oplus\lambda y$. Then by (1) we have
\[ \rho(x,z)+\rho(z,y)=\lambda\rho(x,y)+(1-\lambda)\rho(x,y)=\rho(x,y). \qedhere \]
\end{proof}

All the theorems we will discuss are set in hyperbolic spaces. A metric space $(X,\rho)$ with a coherent system of geodesics is \emph{hyperbolic} if for every $x,y,z\in X$ and every $t\in[0,1]$ we have
\[ \rho\bigl((1-t)x\oplus ty,(1-t)x\oplus tz\bigr)\leq t\rho(y,z). \]
Intuitively, this inequality conveys the idea that triangles in $X$ do not look thicker than Euclidean triangles. It is easy to see that balls are convex sets in hyperbolic spaces. Indeed, for every $x,y,z\in X$ we have
\begin{align*} 
\rho\bigl((1-t)y\oplus tz,x\bigr) &\leq \rho\bigl((1-t)y\oplus tz,(1-t)y\oplus tx\bigr)+\rho\bigl((1-t)y\oplus tx,x\bigr)\leq \\ 
                                  &\leq t\rho(z,x)+(1-t)\rho(y,x)\leq\max\{\rho(y,x),\rho(z,x)\}.
\end{align*}
The following, easy lemma is an immediate consequence of this observation.
\begin{lemma}
\label{lemma:triangle}
Let $(X,\rho)$ be a hyperbolic metric space and let $a,b,c\in X$ and $C>0$ be such that $\rho(a,b)\leq C$, $\rho(b,c)\leq C$ and $\rho(c,a)\leq C$. Then $\rho(x,y)\leq C$ for every $x\in[a,b]$ and every $y\in[b,c]$.
\end{lemma}
\begin{proof}
Since balls are convex, we have $\rho(x,b)\leq C$ and $\rho(x,c)\leq C$. Using the convexity of balls once again, we get $\rho(x,y)\leq C$.
\end{proof}
\noindent Note further that this definition of hyperbolic spaces includes every convex subset of any Banach space. A coherent system of geodesics is given in this case by the usual line segments.

A function $\func{\omega}{[0,+\infty)}{[0,+\infty)}$ is \emph{subadditive} if $\omega(a+b)\leq\omega(a)+\omega(b)$ for every $a,b\in[0,+\infty)$. Given a metric space $(X,\rho)$, recall that a mapping $\func{f}{X}{X}$ is \emph{uniformly continuous} if for every $\epsilon>0$ there is $\delta>0$ such that $\rho(f(x),f(y))<\epsilon$ holds for every $x,y\in X$ with $\rho(x,y)<\delta$. If $(X,\rho)$ is a metric space with a coherent system of geodesics and $\func{f}{X}{X}$ is a uniformly continuous mapping, it is not hard to show that for every $s\in[0,+\infty)$ the values
\[ \omega_f(s)=\sup\bigl\{\rho\bigl(f(x),f(y)\bigr)\,:\,x,y\in X,\,\rho(x,y)\leq s\bigr\} \]
are finite and define a subadditive, nondecreasing function $\func{\omega_f}{[0,+\infty)}{[0,+\infty)}$. Moreover, it is also easy to see that $\omega_f$ is continuous at $0$ and that $\omega_f(0)=0$. Subadditivity and continuity at $0$ imply that $\omega_f$ is continuous everywhere on $[0,+\infty)$. We call $\omega_f$ the \emph{modulus of continuity} of $f$. A family $\mathcal{F}$ of mappings $\func{f}{X}{X}$ is \emph{equicontinuous} if for every $\epsilon>0$ there is $\delta>0$ such that $\rho\bigl(f(x),f(y)\bigr)<\epsilon$ holds whenever $f\in\mathcal{F}$ and $\rho(x,y)<\delta$. Given an equicontinuous family $\mathcal{F}$, it is clear that all its elements are uniformly continuous and that the values $\omega(s)=\sup_{f\in\mathcal{F}}\omega_f(s)$ are finite for every $s\in[0,+\infty)$ and define a subadditive, nondecreasing function $\func{\omega}{[0,+\infty)}{[0,+\infty)}$. It is also easy to see that $\omega(0)=0$ and that $\omega$ is continuous at $0$, which implies that it is continuous everywhere. 

In this note we consider a complete, unbounded hyperbolic space $X$ and a continuous, subadditive, nonzero and nondecreasing function $\func{\omega}{[0,+\infty)}{[0,+\infty)}$ with $\omega(0)=0$ and study the space $\mathcal{C}_\omega(X)$ of all mappings $\func{f}{X}{X}$ with $\omega_f\leq\omega$. $\mathcal{C}_\omega(X)$ can be seen as the largest equicontinuous family of mappings on $X$ whose moduli of continuity have $\omega$ as a common upper bound. The space $\mathcal{C}_\omega(X)$ can be endowed with the topology of uniform convergence on bounded sets. This topology is metrisable by several complete metrics, which make the use of Baire categories meaningful thanks to Baire's theorem. The metric $\func{d}{\mathcal{C}_\omega(X)\times\mathcal{C}_\omega(X)}{\R}$ we use is defined by
\[ d(f,g)=\sum_{n=1}^\infty 2^{-n}\sup_{\rho(x,x_0)\leq n}\Bigl(\min\bigl\{1,\rho\bigl(f(x),g(x)\bigr)\bigr\}\Bigr), \]
where $x_0\in X$ is fixed. The main result will also assume that $\omega$ is concave, i.e.\ $\omega((1-\lambda)a+\lambda b)\geq(1-\lambda)\omega(a)+\lambda\omega(b)$ for every $a,b\in[0,+\infty)$ and every $\lambda\in[0,1]$. The concavity of $\omega$ is needed to make use of the following lemma.
\begin{lemma}
\label{lemma:conc}
Let $\func{\omega}{[0,+\infty)}{[0,+\infty)}$ be a concave function with $\omega(0)=0$. Then $\omega(\lambda s)\leq\lambda\omega(s)$ for every $\lambda\geq 1$ and $s\in[0,+\infty)$.
\end{lemma}
\begin{proof}
We have
\begin{align*}
\omega(\lambda s) &= \lambda\lambda^{-1}\omega(\lambda s)+\lambda(1-\lambda^{-1})\omega(0)\leq \\
                  &\leq \lambda\omega\bigl(\lambda^{-1}\lambda s+(1-\lambda^{-1})0\bigr)=\lambda\omega(s). \qedhere
\end{align*}
\end{proof}
\noindent The concavity of $\omega$ is a stronger requirement than subadditivity. Indeed, it follows by Lemma \ref{lemma:conc} that 
\[ \frac{a}{a+b}\omega(a+b)\leq\omega(a),\quad\frac{b}{a+b}\omega(a+b)\leq\omega(b) \]
for all $a,b\in(0,+\infty)$, which imply $\omega(a+b)\leq\omega(a)+\omega(b)$ for all $a,b\in[0,+\infty)$. Unfortunately, we could not find a way to avoid using this simple lemma in the proof of Theorem \ref{theo:unifb}. This leaves the following question open.
\begin{quest}
Does the conclusion of Theorem \ref{theo:unifb} hold true if $\omega$ is only assumed to be continuous and subadditive, but not concave?
\end{quest}

Finally, we wish to recall some basic facts about Baire categories since, as it has already been mentioned, they play a central role in the statement of the main theorem. Recall that a set $E$ in a topological space $X$ is \emph{nowhere dense} if the closure of $E$ has empty interior in $X$. A set is \emph{meagre} or of the \emph{first category} if it is a countable union of nowhere dense sets. The complement of a meagre set is often said to be \emph{comeagre} or \emph{residual}. Baire's theorem asserts that residual sets are never empty in complete metric spaces. Since a countable intersection of residual sets is still residual, in a complete metric space it is impossible to obtain the empty set by intersecting countably many residual sets. In this sense residual sets are intuitively thought as being large, whereas meagre sets are tiny. Following this intuition, we say that a property is \emph{generic} if it holds at every point of a residual set. Although meagreness is a ubiquitous notion of smallness in complete metric spaces, sometimes one can prove stronger results by using $\sigma$-porosity instead. Following Zaj\'{i}\v{c}ek's definition (see \cite{zajicek_2005}), a set $P$ in a metric space $(X,d)$ is \emph{lower porous} at $x\in X$ if there are $\epsilon_0>0$ and $\alpha>0$ such that for every $\epsilon\in(0,\epsilon_0]$ there is $y$ with $B(y,\alpha\epsilon)\subseteq B(x,\epsilon)\setminus P$. $P$ is said to be lower porous if it is porous at $x$ for every $x\in X$. A set is \emph{$\sigma$-lower porous} if it is a countable union of lower porous sets. Lower porous sets are clearly nowhere dense, hence $\sigma$-lower porous sets are meagre.

\section{Perturbations of uniformly continuous mappings}
\label{sec:main}
The goal of this section is to describe how we can modify uniformly continuous mappings to achieve the main genericity result of the present paper. We start with a lemma.
\begin{lemma}
\label{lemma:omega}
Let $\func{\omega}{[0,+\infty)}{[0,+\infty)}$ be a nondecreasing function such that $\omega(s)\to+\infty$ as $s\to+\infty$. For every $k>0$ and every $\lambda\in(0,1)$, there is $M>0$ such that the inequality $\omega(M+s)\geq k+(1-\lambda)\omega(s)$ holds for every $s\geq 0$.
\end{lemma}
\begin{proof}
We choose $M$ with $\omega(M)\geq\lambda^{-1}k$ and consider two cases.
\begin{enumerate}
\item If $s\geq 0$ is such that $\omega(s)\geq\lambda^{-1}k$, we have
\[ \frac{\omega(M+s)}{\omega(s)}\geq 1\geq\frac{k}{\omega(s)}+1-\lambda, \]
which is what we want. 
\item If $s\geq 0$ is such that $\omega(s)<\lambda^{-1}k$, we have
\[ \omega(M+s)\geq\omega(M)\geq\frac{k}{\lambda}=k+(1-\lambda)\frac{k}{\lambda}>k+(1-\lambda)\omega(s), \]
which is again the conclusion we want. \qedhere
\end{enumerate}
\end{proof}

We also need the following lemma, which appears with an equivalent statement as Lemma 5.3 in \cite{brt_2023}. A previous version for normed spaces had already been discovered by M.\ Dymond and can be found as Lemma 3.1 in \cite{dym_2022}.
\begin{lemma}
\label{lemma:michael}
Let $X$ be a hyperbolic space, let $z_0\in X$ and let $\delta$, $r$ be such that $0<\delta<r$. There exists a map $\func{\Phi}{X}{X}$ which fulfills the following requirements.
\begin{enumerate}
\item $\Phi(x)=z_0$ for every $x\in X$ with $\rho(x,z_0)\leq\delta$.
\item $\Phi(x)=x$ for every $x\in X$ with $\rho(x,z_0)\geq r$.
\item $\Phi$ is a Lipschitz map with Lipschitz constant
	\[ \textup{Lip}\,\Phi\leq\frac{r}{r-\delta}. \]
\end{enumerate} 
\end{lemma}

We are now ready to state and prove the main theorem of this work. 
\begin{theo}
\label{theo:unifb}
Let $(X,\rho)$ be a complete, unbounded hyperbolic metric space and let $\func{\omega}{[0,+\infty)}{[0,+\infty)}$ be a concave, nonzero and nondecreasing function with $\omega(0)=0$. Then there is a meagre set $\mathcal{M}\subset\mathcal{C}_\omega(X)$ such that $\omega_f=\omega$ for every $f\in\mathcal{C}_\omega(X)\setminus\mathcal{M}$.
\end{theo}
\begin{proof}
Let $x_0\in X$ be fixed from now on. Note that the subadditivity of $\omega$ and the requirement that $\omega$ is nonzero and nondecreasing imply that $\omega(s)>0$ for every $s\in(0,+\infty)$. For every $s\in\Q^+$ and every $\mu\in\Q\cap(0,1)$ we consider the set
\[ N(s,\mu)=\{ f\in\mathcal{C}_\omega(X)\,:\,\omega_f(s)\leq(1-\mu)\omega(s) \}. \]
We show that $N(s,\mu)$ is nowhere dense for every $s\in\Q^+$ and $\mu\in\Q\cap(0,1)$. To this end, we fix $s\in\Q^+$, $\mu\in\Q\cap(0,1)$ and a function $f\in\mathcal{C}_\omega(X)$. We also pick an arbitrary $\epsilon>0$ and divide the proof in two steps. \\

\noindent\textbf{Step 1} We find a function $h\in\mathcal{C}_\omega(X)$ such that $d(f,h)<\epsilon$, where $d$ is the distance on $\mathcal{C}_\omega(X)$ defined in Section \ref{sec:setting}, and two points $y_0,z_0\in X$ with $\rho(y_0,z_0)=s$ and $\rho(h(y_0),h(z_0))\geq(1-2^{-1}\mu)\omega(s)$. Clearly, such a function $h$ belongs to $\mathcal{C}_\omega(X)\setminus N(s,\mu)$. The construction of $h$ must be carried out in two different ways, depending on whether $\omega$ is bounded or not. \\

\noindent\textsc{Construction of $h$ with unbounded $\omega$}. Choose a positive integer $p$ with $\sum_{n=p+1}^\infty 2^{-n}<\epsilon/2$. Choose then $t\in(0,1/2)$ with $\omega(p)t<\epsilon/2$. Let $M$ be given by Lemma \ref{lemma:omega} applied to $\omega$, $k=\omega(s)$ and $\lambda=t/2$. Finally, we introduce $R=t^{-1}(2-t)(M+s)$, so that $R>M+s$ and
\[ (1-t)\frac{R}{R-(M+s)}=1-\frac{t}{2}. \] 
We choose $z_0\in X$ with $\rho(z_0,x_0)>p+R$, which exists because $X$ is unbounded. Also set $w_0=(1-t)f(z_0)\oplus tf(x_0)$ and find a point $e_0\in X$ with $\rho(w_0,e_0)=\omega(s)$. Finally, let $\Phi$ be the mapping given by Lemma \ref{lemma:michael} applied to $X$, $z_0$, $\delta=M+s$ and $r=R$. 

We define a function $\func{g}{X}{X}$ by 
\[ g(x)=(1-t)f\bigl(\Phi(x)\bigr)\oplus tf(x_0) \]
and note that $g(x)=w_0$ for every $x\in X$ such that $\rho(x,z_0)\leq M+s$. We want to show $\omega_g\leq(1-t/2)\omega$. Picking $x,y\in X$ we have
\begin{align*}
\rho\bigl(g(x),g(y)\bigr) &= \rho\bigl((1-t)f\circ\Phi(x)\oplus tf(x_0),(1-t)f\circ\Phi(y)\oplus tf(x_0)\bigr)\leq \\
 			  &\leq (1-t)\rho\bigl(f\circ\Phi(x),f\circ\Phi(y)\bigr)\leq \\
		          &\leq (1-t)\omega\Bigl(\rho\bigr(\Phi(x),\Phi(y)\bigr)\Bigr)\leq \\
		          &\leq (1-t)\omega\biggl(\frac{R}{R-(M+s)}\rho(x,y)\biggr)\leq \\
		          &\leq (1-t)\frac{R}{R-(M+s)}\omega\bigl(\rho(x,y)\bigr)= \\
		          &= \biggl(1-\frac{t}{2}\biggr)\omega\bigl(\rho(x,y)\bigr),
\end{align*}
where the first inequality follows from $X$ being hyperbolic and the last two steps in the computation follow from Lemma \ref{lemma:conc} and the definition of $R$ respectively. 

We define the function $\func{h}{X}{X}$ we are looking for in the following way
\[ h(x)=\left\{\begin{array}{ll}
		g(x) & \text{if }\rho(x,z_0)\geq s \\
		\Bigl(1-\frac{\omega(s-\rho(x,z_0))}{\omega(s)}\Bigr)w_0\oplus\frac{\omega(s-\rho(x,z_0))}{\omega(s)}e_0 & \text{if }\rho(x,z_0)<s
		\end{array}\right. \]
We have to check that $\omega_h\leq\omega$. We pick $x,y\in X$ and consider four cases. The other cases are dealt with symmetrically by swapping the roles of $x$ and $y$. \\

\noindent\emph{Case 1.} If $\rho(x,z_0)\geq s$ and $\rho(y,z_0)\geq s$, then $h(x)=g(x)$ and $h(y)=g(y)$, therefore there is nothing to show. \\

\noindent\emph{Case 2.} If $\rho(x,z_0)<s$ and $s\leq\rho(y,z_0)\leq M+s$, we have $h(y)=w_0$, hence
\[ \rho\bigl(h(x),h(y)\bigr)=\omega\bigl(s-\rho(x,z_0)\bigr)\leq\omega\bigl(\rho(y,z_0)-\rho(x,z_0)\bigr)\leq\omega\bigl(\rho(x,y)\bigr) \]

\noindent\emph{Case 3.} If $\rho(x,z_0)\leq\rho(y,z_0)<s$, then
\begin{align*}
\rho\bigl(h(x),h(y)\bigr) &= \omega\bigl(s-\rho(x,z_0)\bigr)-\omega\bigl(s-\rho(y,z_0)\bigr)\leq \\
 		          &\leq\omega\bigl(\rho(y,z_0)-\rho(x,z_0)\bigr)\leq \\
 		          &\leq\omega\bigl(\rho(x,y)\bigr).
\end{align*}
\noindent\emph{Case 4.} If $\rho(x,z_0)<s$ and $\rho(y,z_0)>M+s$, by the continuity of the distance to $z_0$ and the geodesic from $x$ to $y$ there must be two points $a$ and $b$ on the geodesic segment $[x,y]$ such that $\rho(a,z_0)=s$ and $\rho(b,z_0)=M+s$. Note that $h(a)=h(b)=w_0$. On one hand we have 
\begin{align*}
\omega\bigl(\rho(x,y)\bigr) &= \omega\bigl(\rho(x,a)+\rho(a,b)+\rho(b,y)\bigr)\geq \\
			    &= \omega\bigl(\rho(a,b)+\rho(b,y)\bigr)\geq \\
		            &\geq \omega\bigl(\rho(b,z_0)-\rho(a,z_0)+\rho(b,y)\bigr)= \\
		            &= \omega\bigl(M+s-s+\rho(b,y)\bigr)=\omega\bigl(M+\rho(b,y)\bigr),
\end{align*}
where the first equality follows from Lemma \ref{lemma:segment}. On the other hand we have
\begin{align*}
\rho\bigl(h(x),h(y)\bigr) &\leq \rho\bigl(h(x),h(a)\bigr)+\rho\bigl(h(a),h(b)\bigr)+\rho\bigl(h(b),h(y)\bigr)= \\
		          &= \omega\bigl(s-\rho(x,z_0)\bigr)+\rho\bigl(g(b),g(y)\bigr)\leq \\
		          &\leq \omega(s)+\biggl(1-\frac{t}{2}\biggr)\omega\bigl(\rho(b,y)\bigr).
\end{align*}
Hence, we are done if 
\[ \omega\bigl(M+\rho(b,y)\bigr)\geq \omega(s)+\biggl(1-\frac{t}{2}\biggr)\omega\bigl(\rho(b,y)\bigr), \]
but this is true by the choice of $M$. \\

At this point, pick any $y_0$ with $\rho(y_0,z_0)=s$, which exists because $X$ is unbounded and connected. Now we have 
\[ \rho\bigl(h(z_0),h(y_0)\bigr)=\rho(e_0,w_0)=\omega(s)\geq\biggl(1-\frac{\mu}{2}\biggr)\omega(s), \]
as wished. \\

\noindent\textsc{Construction of $h$ with bounded $\omega$.} Set $\Omega=\sup_{u>0}\omega(u)$. As before, choose a positive integer $p$ with $\sum_{n=p+1}^\infty 2^{-n}<\epsilon/2$. Choose then $t>0$ with
\[ t\leq\min\biggl\{\frac{\mu}{4},\frac{\epsilon}{2\Omega}\biggr\}. \]
We first define $\func{g_1}{X}{X}$ by $g_1(x)=(1-t)f(x)\oplus tf(x_0)$ for every $x\in X$ and note that $\omega_{g_1}\leq(1-t)\omega_f$. Indeed, since $X$ is hyperbolic we have
\begin{align*}
\rho\bigl(g_1(x),g_1(y)\bigr) &= \rho\bigl((1-t)f(x)\oplus tf(x_0),(1-t)f(y)\oplus tf(x_0)\bigr)\leq \\
			  &\leq (1-t)\rho\bigl(f(x),f(y)\bigr).
\end{align*} 
Let $s'$ be such that $\omega(s')\geq(1-t)\Omega$. Set $M=s+3s'$ and choose $R\geq t^{-1}M$, so that 
\[ (1-t)\frac{R}{R-M}\leq 1. \]
Choose a point $z_0\in X$ such that $\rho(z_0,x_0)>p+R$ and set $w_0=g_1(z_0)$. We now need to consider the function
\[ \func{\phi}{X}{[0,+\infty)},\quad \phi(x)=\sup_{y\in g_1(X)}\rho(x,y). \]
Note that $\phi(x)\leq(1-t)\Omega$ for every $x\in g_1(X)$. Since $\phi$ is continuous and $X$ is connected and unbounded, we have $[(1-t)\Omega,+\infty)\subseteq\phi(X)$, hence there must be $e_0\in X$ with $\phi(e_0)=(1-t)\Omega$. Find $w_1\in g_1(X)$ with $\rho(w_1,e_0)>(1-2t)\Omega$. Further, observe that $\rho(w_0,w_1)\leq(1-t)\Omega$, as both points lie in $g_1(X)$. Let $\func{\Phi}{X}{X}$ be the mapping given by Lemma \ref{lemma:michael} applied to $z_0$, $\delta=M$ and $r=R$. Define $\func{g_2}{X}{X}$ by $g_2=g_1\circ\Phi$ and note that $g_2(x)=w_0$ for every $x\in X$ with $\rho(x,z_0)\leq M$. Moreover, again by Lemma \ref{lemma:conc} we have 
\begin{align*}
\rho\bigl(g_2(x),g_2(y)\bigr) &\leq \omega_{g_1}\bigl(\rho(\Phi(x),\Phi(y))\bigr)\leq \\
			      &\leq (1-t)\omega\biggl(\frac{R}{R-M}\rho(x,y)\biggl)\leq \\
			      &\leq (1-t)\frac{R}{R-M}\omega\bigl(\rho(x,y)\bigr)\leq \\
			      &\leq \omega\bigl(\rho(x,y)\bigr)
\end{align*}
for every $x,y\in X$, therefore $g_2\in\mathcal{C}_\omega(X)$. Finally, define the point
\[ w_2=\biggl(1-\frac{\omega(s)}{\Omega}\biggr)w_1\oplus\frac{\omega(s)}{\Omega}e_0. \]
We are now ready to define the function $\func{h}{X}{X}$ we are looking for.
\begin{itemize}
\item[(i)] If $\rho(x,z_0)>s+2s'$, we set $h(x)=g_2(x)$.
\item[(ii)] If $s+s'<\rho(x,z_0)\leq s+2s'$, we set
\[ h(x)=\biggl(1-\frac{\omega(s+2s'-\rho(x,z_0))}{\omega(s')}\biggr)w_0\oplus\frac{\omega(s+2s'-\rho(x,z_0))}{\omega(s')}w_1. \]
\item[(iii)] If $s<\rho(x,z_0)\leq s+s'$, we set $h(x)=w_1$.
\item[(iv)] If $\rho(x,z_0)\leq s$, we set 
\[ h(x)=\biggl(1-\frac{\omega(s-\rho(x,z_0))}{\omega(s)}\biggr)w_1\oplus\frac{\omega(s-\rho(x,z_0))}{\omega(s)}w_2. \]
\end{itemize}
Let us check that $h\in\mathcal{C}_\omega(X)$. We pick $x,y\in X$ and consider several cases. The remaining cases are checked symmetrically by swapping $x$ and $y$.\\

\noindent\emph{Case 1.} If $\rho(x,z_0)\geq s+2s'$ and $\rho(y,z_0)\geq s+2s'$, then $h(x)=g_2(x)$ and $h(y)=g_2(y)$, hence there is nothing to show. \\

\noindent\emph{Case 2.} If $s+s'<\rho(x,z_0)\leq\rho(y,z_0)\leq s+2s'$, we have
\begin{align*}
\rho\bigl(h(x),h(y)\bigr) &= \biggl(\frac{\omega(s+2s'-\rho(x,z_0))}{\omega(s')}-\frac{\omega(s+2s'-\rho(y,z_0))}{\omega(s')}\biggr)\rho(w_0,w_1)\leq \\
			  &\leq \frac{\rho(w_0,w_1)}{\omega(s')}\omega\bigl(\rho(y,z_0)-\rho(x,z_0)\bigr)\leq \\
			  &\leq \frac{(1-t)\Omega}{\omega(s')}\omega\bigl(\rho(x,y)\bigr)\leq \\
			  &\leq \omega\bigl(\rho(x,y)\bigr).
\end{align*}
\vspace{1pt}

\noindent\emph{Case 3.} If $s+s'<\rho(x,z_0)\leq s+2s'$ and $s+2s'<\rho(y,z_0)\leq s+3s'$, then $h(y)=w_0$, therefore
\begin{align*}
\rho\bigl(h(x),h(y)\bigr) &= \frac{\omega(s+2s'-\rho(x,z_0))}{\omega(s')}\rho(w_0,w_1)\leq \\
			  &\leq \frac{\rho(w_0,w_1)}{\omega(s')}\omega\bigl(\rho(y,z_0)-\rho(x,z_0)\bigr)\leq \\
			  &\leq \frac{(1-t)\Omega}{\omega(s')}\omega\bigl(\rho(x,y)\bigr)\leq \\
			  &\leq \omega\bigl(\rho(x,y)\bigr).
\end{align*}
\vspace{1pt}

\noindent\emph{Case 4.} If $s+s'<\rho(x,z_0)\leq s+2s'$ and $\rho(y,z_0)>s+3s'$, then $h(y)\in g_1(X)$, which implies $\rho(w_0,h(y))\leq(1-t)\Omega$ and $\rho(w_1,h(y))\leq(1-t)\Omega$. Since balls are convex and $h(x)\in[w_0,w_1]$, we have
\[ \rho\bigl(h(x),h(y)\bigr)\leq(1-t)\Omega\leq\omega(s')\leq\omega\bigl(\rho(y,z_0)-\rho(x,z_0)\bigr)\leq\omega\bigl(\rho(x,y)\bigr). \]
\vspace{1pt}

\noindent\emph{Case 5.} If $s<\rho(x,z_0)\leq s+s'$ and $s<\rho(y,z_0)\leq s+s'$, then $h(x)=h(y)=w_1$ and there is nothing to show. \\

\noindent\emph{Case 6.} If $s<\rho(x,z_0)\leq s+s'$ and $s+s'<\rho(y,z_0)\leq s+2s'$, then $h(x)=w_1$ and
\begin{align*}
\rho\bigl(h(x),h(y)\bigr) &= \biggl(1-\frac{\omega(s+2s'-\rho(y,z_0))}{\omega(s')}\biggr)\rho(w_0,w_1)= \\
			  &= \frac{\rho(w_0,w_1)}{\omega(s')}\Bigl(\omega(s')-\omega\bigl(s+2s'-\rho(y,z_0)\bigr)\Bigr)\leq \\
			  &\leq \frac{(1-t)\Omega}{\omega(s')}\omega\bigl(\rho(y,z_0)-(s+s')\bigr)\leq \\
			  &\leq \omega\bigl(\rho(y,z_0)-\rho(x,z_0)\bigr)\leq \\
			  &\leq \omega\bigl(\rho(x,y)\bigr).
\end{align*}
\vspace{1pt}

\noindent\emph{Case 7.} If $\rho(x,z_0)\leq s+s'$ and $\rho(y,z_0)>s+2s'$, we have $h(y)\in g_1(X)$, hence $\rho(e_0,h(y))\leq(1-t)\Omega$ and $\rho(w_1,h(y))\leq(1-t)\Omega$. Since balls are convex and $h(x)\in[w_1,e_0]$, we have
\[ \rho\bigl(h(x),h(y)\bigr)\leq(1-t)\Omega\leq\omega(s')\leq\omega\bigl(\rho(y,z_0)-\rho(x,z_0)\bigr)\leq\omega\bigl(\rho(x,y)\bigr). \] 
\vspace{1pt}

\noindent\emph{Case 8.} If $\rho(x,z_0)\leq\rho(y,z_0)\leq s$, then
\begin{align*}
\rho\bigl(h(x),h(y)\bigr) &= \biggl(\frac{\omega(s-\rho(x,z_0))}{\omega(s)}-\frac{\omega(s-\rho(y,z_0))}{\omega(s)}\biggr)\rho(w_1,w_2)\leq \\
			  &\leq \frac{\rho(w_1,w_2)}{\omega(s)}\omega\bigl(\rho(y,z_0)-\rho(x,z_0)\bigr)\leq \\
			  &\leq \frac{\omega(s)}{\Omega}\cdot\frac{\rho(w_1,e_0)}{\omega(s)}\omega\bigl(\rho(x,y)\bigr)\leq \\
			  &\leq \frac{\rho(w_1,e_0)}{\Omega}\omega\bigl(\rho(x,y)\bigr)\leq \\
			  &\leq \omega\bigl(\rho(x,y)\bigr).
\end{align*}
\vspace{1pt}

\noindent\emph{Case 9.} If $\rho(x,z_0)\leq s$ and $s<\rho(y,z_0)\leq s+s'$, then $h(y)=w_1$ and
\begin{align*}
\rho\bigl(h(x),h(y)\bigr) &= \frac{\omega(s-\rho(x,z_0))}{\omega(s)}\rho(w_1,w_2)\leq \\
			  &\leq \frac{\rho(w_1,e_0)}{\Omega}\omega\bigl(\rho(y,z_0)-\rho(x,z_0)\bigr)\leq \\
			  &\leq \omega\bigl(\rho(x,y)\bigr).
\end{align*}
\vspace{1pt}

\noindent\emph{Case 10.} If $\rho(x,z_0)\leq s$ and $s+s'<\rho(y,z_0)\leq s+2s'$, then $h(x)\in[w_1,e_0]$ and $h(y)\in[w_0,w_1]$, hence it follows by Lemma \ref{lemma:triangle} that
\[ \rho\bigl(h(x),h(y)\bigr)\leq(1-t)\Omega\leq\omega(s')\leq\omega\bigl(\rho(y,z_0)-\rho(x,z_0)\bigr)\leq\omega\bigl(\rho(x,y)\bigr). \]
\vspace{1pt}

Pick now any $y_0$ with $\rho(y_0,z_0)=s$. We have
\begin{align*} 
\rho\bigl(h(y_0),h(z_0)\bigr) &= \rho(w_1,w_2)=\frac{\omega(s)}{\Omega}\rho(w_1,e_0)> \\
			      &> \frac{\omega(s)}{\Omega}(1-2t)\Omega=(1-2t)\omega(s)\geq\biggl(1-\frac{\mu}{2}\biggr)\omega(s),
\end{align*} 
hence we found $h$, $z_0$ and $y_0$ as wished. \\

Now that $h$ has been constructed, it remains to check that $d(f,h)<\epsilon$. To see this, note that $h(x)=(1-t)f(x)\oplus tf(x_0)$ whenever $x\in X$ is such that $\rho(x,x_0)\leq p$, therefore
\begin{align*}
d(f,h) &\leq \sum_{n=1}^p2^{-n}\sup_{\rho(x,x_0)\leq n}\Bigl(\min\bigl\{1,\rho\bigl(f(x),g_1(x)\bigr)\bigr\}\Bigr)+\sum_{n=p+1}^\infty 2^{-n}< \\
       &< \sum_{n=1}^p2^{-n}t\rho\bigl(f(x),f(x_0)\bigr)+\frac{\epsilon}{2}< \\
       &< t\,\Omega\cdot\sum_{n=1}^\infty 2^{-n}+\frac{\epsilon}{2}\leq \\
       &\leq \frac{\epsilon}{2}+\frac{\epsilon}{2}=\epsilon. 
\end{align*}
\vspace{1pt}

\noindent\textbf{Step 2.} We want to show that there is $\eta>0$ such that $h'\notin N(s,\mu)$ whenever $h'\in\mathcal{C}_\omega(X)$ and $d(h,h')<\eta$, where $h$ is the mapping we constructed in Step 1. Find a positive integer $q$ such that $y_0,z_0\in B(x_0,q)$. Set
\[ \eta=2^{-q}\min\biggl\{1,\frac{\omega(s)\mu}{4}\biggr\}. \]
If $h'\in\mathcal{C}_\omega(X)$ and $d(h,h')<\eta$ we have
\begin{align*} 
\eta &> \sum_{n=1}^\infty 2^{-n}\sup_{\rho(x,x_0)\leq n}\Bigl(\min\bigl\{1,\rho\bigl(h(x),h'(x)\bigr)\bigr\}\Bigr)\geq \\
     &\geq 2^{-q}\min\bigl\{1,\rho\bigl(h(y_0),h'(y_0)\bigr)\bigr\}.
\end{align*}
This can happen only if $\rho(h(y_0),h'(y_0))<4^{-1}\omega(s)\mu$. Similarly, one shows that $\rho(h(z_0),h'(z_0))<4^{-1}\omega(s)\mu$. It follows that
\begin{align*}
\omega_{h'}(s) &\geq \rho\bigl(h'(y_0),h'(z_0)\bigr)\geq \\
	       &\geq \rho\bigl(h(y_0),h(z_0)\bigr)-\rho\bigl(h(z_0),h'(z_0)\bigr)-\rho\bigl(h(y_0),h'(y_0)\bigr)> \\
	       &> \biggl(1-\frac{\mu}{2}\biggr)\omega(s)-\frac{\omega(s)\mu}{4}-\frac{\omega(s)\mu}{4}=(1-\mu)\omega(s).
\end{align*}
That is, $h'\notin N(s,\mu)$, as we desire. \\

It follows from Step 1 and Step 2 that $\mathcal{C}_\omega(X)\setminus N(s,\mu)$ contains a dense, open set, which is equivalent to $N(s,\mu)$ being nowhere dense, as we wanted. Now we define
\[ \mathcal{M}=\bigcup_{s\in\Q^+}\bigcup_{\mu\in\Q\cap(0,1)}N(s,\mu), \]
which is meagre because it is a countable union of nowhere dense sets. Clearly, if $f\in\mathcal{C}_\omega(X)\setminus\mathcal{M}$ we have $\omega_f(s)=\omega(s)$ for every $s\in\Q^+$. By the continuity of $\omega_f$ and $\omega$, we also have $\omega_f(s)=\omega(s)$ for every $s\in[0,+\infty)$, hence the proof is complete.
\end{proof}

\section{Spaces of bounded mappings}
Let $X$ be a complete hyperbolic space and let $\func{\omega}{[0,+\infty)}{[0,+\infty)}$ be continuous, subadditive, nondecreasing and such that $\omega(0)=0$. In addition to the space $\mathcal{C}_\omega(X)$ we have already discussed, we can also consider
\[ \mathcal{C}_\omega^\textup{b}(X)=\Bigl\{ f\in\mathcal{C}_\omega(X)\,:\,\sup_{x,y\in X}\rho\bigl(f(x),f(y)\bigr)<+\infty\Bigr\} \]
and endow this space with the topology of uniform convergence given by the metric 
\[ \func{d_\infty}{\mathcal{C}_\omega^\textup{b}(X)\times\mathcal{C}_\omega^\textup{b}(X)}{\R},\quad d_\infty(f,g)=\sup_{x\in X}\rho\bigl(f(x),g(x)\bigr). \]
If $\omega$ is bounded, then $\mathcal{C}_\omega(X)=\mathcal{C}_\omega^\textup{b}(X)$ as sets, but the topology of uniform convergence and the topology of uniform convergence on bounded sets are different. The meagreness result of Section \ref{sec:main} also changes drastically. The proof of the next theorem follows the same ideas one can find in \cite{rz_2001}, Theorem 2.2, and does not need to assume that $\omega$ is concave.
\begin{theo}
\label{theo:uniftop}
Let $(X,\rho)$ be a complete, unbounded hyperbolic metric space and let $\func{\omega}{[0,+\infty)}{[0,+\infty)}$ be a continuous, subadditive, nonzero and nondecreasing function with $\omega(0)=0$. Then there is a $\sigma$-lower porous set $\mathcal{M}\subset\mathcal{C}^\textup{b}_\omega(X)$ such that for every $f\in\mathcal{C}_\omega^\textup{b}(X)\setminus\mathcal{M}$ one has $\omega_f(s)<\omega(s)$ at every point $s\in(0,+\infty)$.
\end{theo}
\begin{proof}
Let $x_0\in X$ be fixed. For every $s\in\Q^+$ we consider the set $R(s)$ of mappings $f\in\mathcal{C}_\omega^\textup{b}(X)$ which satisfy the following property: there is $\mu\in(0,1)$ such that $\rho(f(x),f(y))\leq\mu\,\omega(\rho(x,y))$ for every $x,y\in X$ with $\rho(x,y)\geq s$. The goal is to show that $\mathcal{C}_\omega^\textup{b}(X)\setminus R(s)$ is lower porous for every $s\in\Q^+$. We now fix $s\in\Q^+$. Given $f\in\mathcal{C}^\textup{b}_\omega(X)$, we want $\epsilon_0>0$ and $\alpha>0$ such that, for every $\epsilon\in(0,\epsilon_0]$, there is $g\in\mathcal{C}^\textup{b}_\omega(X)$ with $B(g,\alpha\epsilon)\subseteq B(f,\epsilon)\cap R(s)$. Define $\Omega=\sup_{t>0}\omega_f(t)$. \\

\noindent\textbf{Case 1.} If $\Omega>0$, set $\epsilon_0=2\Omega$. Given $\epsilon\in(0,\epsilon_0]$, define $\gamma=2^{-1}\Omega^{-1}\epsilon$ and set $\func{g}{X}{X}$ to be
\[ g(x)=(1-\gamma)f(x)\oplus\gamma f(x_0). \]
Let us check that $d_\infty(f,g)\leq\epsilon/2$. We have
\[ d_\infty(f,g)=\sup_{x\in X}\gamma\rho\bigl(f(x),f(x_0)\bigr)\leq\gamma\,\Omega=\frac{\epsilon}{2}. \]
Moreover, since $X$ is hyperbolic we have $\omega_g(t)\leq(1-\gamma)\omega_f(t)\leq(1-\gamma)\omega(t)$ for every $t>0$. Next, set 
\[ \alpha=\min\biggl\{\frac{1}{2},\frac{\omega(s)}{8\,\Omega}\biggr\} \] 
and let $h\in\mathcal{C}_\omega^\textup{b}(X)$ be such that $d_\infty(g,h)<\alpha\epsilon$. Clearly, $d_\infty(f,h)<\epsilon$. If $x,y\in X$ are such that $\rho(x,y)\geq s$, we have
\begin{align*}
\rho\bigl(h(x),h(y)\bigr) &\leq \rho\bigl(h(x),g(x)\bigr)+\rho\bigl(g(x),g(y)\bigr)+\rho\bigl(g(y),h(y)\bigr)< \\
	                  &< 2\alpha\epsilon+(1-\gamma)\omega\bigl(\rho(x,y)\bigr)\leq \\
	                  &\leq \frac{\gamma}{2}\omega(s)+(1-\gamma)\omega\bigl(\rho(x,y)\bigr)\leq \\
	                  &\leq \biggl(1-\frac{\gamma}{2}\biggr)\omega\bigl(\rho(x,y)\bigr),
\end{align*}
which implies $h\in R(s)$, hence we have $B(g,\alpha\epsilon)\subseteq B(f,\epsilon)\cap R(s)$, as wished.\\

\noindent\textbf{Case 2.} If $\Omega=0$, then $f$ is constant. Set $\epsilon_0=4^{-1}\omega(s)$ and, for every $\epsilon\in(0,\epsilon_0]$, set $g=f$. Let $\alpha=1$. If $d_\infty(h,g)<\alpha\epsilon=\epsilon$ and $x,y\in X$ are such that $\rho(x,y)\geq s$ we have
\[ \rho\bigl(h(x),h(y)\bigr)\leq\rho\bigl(h(x),g(x)\bigr)+\bigl(g(y),h(y)\bigr)\leq 2\epsilon\leq 2\epsilon_0=\frac{\omega(s)}{2}\leq\frac{1}{2}\omega\bigl(\rho(x,y)\bigr). \]
Hence we have once again $B(g,\alpha\epsilon)\subseteq B(f,\epsilon)\cap R(s)$, as desired. \\

Consider now
\[ \mathcal{R}=\bigcap_{s\in\Q^+} R(s),\quad\mathcal{M}=\mathcal{C}^\textup{b}_\omega(X)\setminus\mathcal{R}, \]
and notice that $\mathcal{M}$ is $\sigma$-lower porous. For $f\in\mathcal{R}$ define $\func{\phi_f}{(0,+\infty)}{\R}$ by
\[ \phi_f(s)=\inf\bigl\{\mu\,:\,\rho\bigl(f(x),f(y)\bigr)\leq\mu\,\omega\bigl(\rho(x,y)\bigr)\text{ for all }x,y\in X\text{ with }\rho(x,y)\geq s\bigr\}. \]
By the choice of $f$ we have $\phi_f(s)<1$ for all $s\in\Q^+$. Since $\phi_f$ is clearly nonincreasing, we also have $\phi_f(s)<1$ for every $s\in(0,+\infty)$. Finally, note that 
\[ \rho\bigl(f(x),f(y)\bigr)\leq\phi_f\bigl(\rho(x,y)\bigr)\omega\bigl(\rho(x,y)\bigr) \] 
for every $x,y\in X$ with $x\ne y$, which implies $\omega_f(s)\leq\phi_f(s)\omega(s)<\omega(s)$ for every $s\in(0,+\infty)$, thus the proof is complete.
\end{proof}

\section{A remark on pointwise convergence}
\label{sec:pointwise}
We consider again $\mathcal{C}_\omega(X)$, where $X$ is a complete, unbounded hyperbolic metric space and $\func{\omega}{[0,+\infty)}{[0,+\infty)}$ is a concave, nonzero, nondecreasing function with $\omega(0)=0$. In Section \ref{sec:main} we endowed $\mathcal{C}_\omega(X)$ with the topology of uniform convergence on bounded sets. However, if $X$ is separable, the topology of pointwise convergence can also be metrised by a complete metric. This is shown in \cite{brt_2023} for nonexpansive mappings, but the argument can easily be extended to uniformly continuous mappings. We briefly present the proof.
\begin{prop}
\label{prop:point}
Let $(X,\rho)$ be a complete, separable metric space and let the set $\Theta={\{\theta_n\}}_{n=1}^\infty$ be dense in $X$. Define $\func{d_\Theta}{\mathcal{C}_\omega(X)\times\mathcal{C}_\omega(X)}{\R}$ by 
\[ d_\Theta(f,g)=\sum_{n=1}^\infty 2^{-n}\min\bigl\{1,\rho\bigl(f(\theta_n),g(\theta_n)\bigr)\bigr\}. \]
Then $d_\Theta$ is a complete metric on $\mathcal{C}_\omega(X)$. Moreover, a sequence ${(f_k)}_{k=1}^\infty$ in $\mathcal{C}_\omega(X)$ converges to $f\in\mathcal{C}_\omega(X)$ if and only if $f_k(x)\to f(x)$ for every $x\in X$ as $k\to\infty$.
\end{prop}
\begin{proof}
Let us show first that $d_\Theta$ is complete. Consider a Cauchy sequence ${(f_k)}_{k=1}^\infty$. It is not hard to see that ${(f_k(\theta_n))}_{k=1}^\infty$ is a Cauchy sequence in $X$ for every $n$, hence it converges to some point $f(\theta_n)$ by the completeness of $X$. This defines a mapping $\func{f}{\Theta}{X}$. One easily checks that $f$ is uniformly continuous on $\Theta$ and that $\omega_f\leq\omega$. Hence, using again the completeness of $X$, we can find a unique extension $f\in\mathcal{C}_\omega(X)$. Since $f_k(\theta_n)$ converges to $f(\theta_n)$ for every $n$, it is straightforward that $d_\Theta(f_k,f)\to 0$ as $k\to\infty$.

If a sequence ${(f_k)}_{k=1}^\infty$ in $\mathcal{C}_\omega(X)$ converges to $f\in\mathcal{C}_\omega(X)$ pointwise, then $f_k(\theta_n)$ converges to $f(\theta_n)$ for every $n$, which implies as above that $d_\Theta(f_k,f)\to 0$ as $k\to\infty$. Conversely, if $d_\Theta(f_k,f)\to 0$ as $k\to\infty$, then $f_k(\theta_n)$ converges to $f(\theta_n)$ for every $n$. Pick an arbitrary $x\in X$ and $\epsilon>0$. Find $n$ such that $\omega(\rho(x,\theta_n))<\epsilon/3$. Find $k$ such that $\rho(f_k(\theta_n),f(\theta_n))<\epsilon/3$. Then we have
\begin{align*}
\rho\bigl(f_k(x),f(x)\bigr) &\leq \rho\bigl(f_k(x),f_k(\theta_n)\bigr)+\rho\bigl(f_k(\theta_n),f(\theta_n)\bigr)+\rho\bigl(f(\theta_n),f(x)\bigr)< \\
			    &< 2\,\omega\bigl(\rho(x,\theta_n)\bigr)+\frac{\epsilon}{3}<\epsilon.
\end{align*}
Since $x\in X$ and $\epsilon>0$ are arbitrary, this proves the claim.
\end{proof}
It turns out that it is possible to fully recover the conclusion of Theorem \ref{theo:unifb} also in the topology of pointwise convergence. We only sketch the proof as it does not exhibit any substantial differences.
\begin{theo}
Let $(X,\rho)$ be a complete, separable, unbounded hyperbolic metric space and let $\func{\omega}{[0,+\infty)}{[0,+\infty)}$ be a concave, nonzero and nondecreasing function with $\omega(0)=0$. If $\mathcal{C}_\omega(X)$ is endowed with the topology of pointwise convergence, then there is a meagre set $\mathcal{M}\subset\mathcal{C}_\omega(X)$ such that $\omega_f=\omega$ for every $f\in\mathcal{C}_\omega(X)\setminus\mathcal{M}$.
\end{theo} 
\begin{proof}
We fix $x_0\in X$ and a dense set $\Theta={\{\theta_n\}}_{n=1}^\infty\subseteq X$ and define the sets $N(\mu,s)$ as in the proof of Theorem \ref{theo:unifb} for every $s\in\Q^+$ and every $\mu\in\Q\cap(0,1)$. The goal is to show that $N(s,\mu)$ is nowhere dense for every $s\in\Q^+$ and $\mu\in\Q\cap(0,1)$. We fix $s\in\Q^+$, $\mu\in\Q\cap(0,1)$ and a function $f\in\mathcal{C}_\omega(X)$. We pick an arbitrary $\epsilon>0$ and divide again the proof in two steps. \\

\noindent\textbf{Step 1.} We want to construct a function $h\in\mathcal{C}_\omega(X)$ such that $d_\Theta(f,h)<\epsilon$, where $d_\Theta$ is the distance defined in the statement of Proposition \ref{prop:point}, and two points $y_0,z_0\in X$ with $\rho(y_0,z_0)=s$ and $\rho(h(y_0),h(z_0))\geq(1-2^{-1}\mu)\omega(s)$. This time we choose $p>0$ such that
\[ \sum_{\rho(\theta_n,x_0)\geq p}2^{-n}<\frac{\epsilon}{2}. \]
The determination of the points $y_0$ and $z_0$ and the construction of $h$ follow now in the exact same way as in the proof of Theorem \ref{theo:unifb} both in the case where $\omega$ is unbounded and in the case where $\omega$ is bounded. \\

\noindent\textbf{Step 2.} We want to exhibit an open neighbourhood $U$ of $h$ such that $h'\notin N(\mu,s)$ for every $h'\in U$. Consider
\begin{align*}
V_1 &= \bigl\{ h'\in\mathcal{C}_\omega(X)\,:\,\rho\bigl(h'(y_0),h(y_0)\bigr)<\omega(s)\mu/4 \bigr\}, \\
V_2 &= \bigl\{ h'\in\mathcal{C}_\omega(X)\,:\,\rho\bigl(h'(z_0),h(z_0)\bigr)<\omega(s)\mu/4 \bigr\}.
\end{align*}
Then $U=V_1\cap V_2$ is the desired neighbourhood. \\

By Step 1 and Step 2, $N(\mu,s)$ is meagre. We define $\mathcal{M}$ as in the proof of Theorem \ref{theo:unifb} and show in the same way that $\mathcal{M}$ is the set we are looking for.
\end{proof}

We conclude by mentioning that the most interesting questions concerning the topology of pointwise convergence are related to nonexpansive mappings and fixed point theory. In particular, the following problem seems to be open.
\begin{quest}
Let $X$ be a complete, separable hyperbolic space and let $\mathcal{N}(X)$ be the space of nonexpansive self-mappings on $X$, endowed with the topology of pointwise convergence. Does the generic element $f\in\mathcal{N}(X)$ have a unique fixed point? If the answer is positive, does the sequence of iterates ${(f^n(x))}_{n=1}^\infty$ converge to the fixed point of $f$ for every $x\in X$? 
\end{quest}

\section*{Acknowledgements}
The author wishes to thank Professor Christian Bargetz for the several, fruitful discussions about this topic and Professor Simeon Reich and the referee for their careful reading of the manuscript and their valuable comments. The author's work has been supported by the Austrian Science Fund (FWF): P 32523-N.

\printbibliography

\end{document}